\documentclass[a4paper,12pt]{article}
\usepackage[T1]{fontenc}
\usepackage{mathptmx}
\usepackage[scaled=0.96]{helvet}
\usepackage{amsmath,amssymb}        
\usepackage[dvips]{graphicx}
\usepackage{psfrag}
\usepackage{url}

\newcommand{\eqnsection}{
   \renewcommand{\theequation}{{\thesection.\arabic{equation}}}
   \makeatletter
   \csname @addtoreset\endcsname{equation}{section}
   \makeatother}
\eqnsection

\def\RR{\hbox{I\kern-.2em\hbox{R}}}
\title{A Non-Iterative Transformation Method for a Class of
Free Boundary Value Problems Governed by ODEs}
\author{Riccardo Fazio\footnote{Corresponding author: e-mail: rfazio@unime.it \ \ \
home-page: http://mat521.unime.it/fazio}\quad and\quad Salvatore Iacono\\[1ex] 
Department of Mathematics and Computer Science \\
University of Messina \\
Viale F. Stagno D'Alcontres 31, 98166 Messina, Italy}
\date{\today}
\begin{document}
\maketitle
\begin{abstract}
The aim of this work is to point out that the class of free boundary problems governed by second order autonomous ordinary differential equations can be transformed to initial value problems.
Interest in the numerical solution of free boundary problems arises because these are always nonlinear problems. 
The theoretical content of this paper is original: results already available in literature are related to the invariance properties of scaling or spiral groups of point transformations but here we show how it is also possible to use e invariance properties of a translation group.
We test the proposed algorithm by solving three problems: a problem describing a rope configuration against an obstacle, a dynamical problem with a nonlinear force, and a problem related to the optimal length estimate for tubular flow reactors. 
\end{abstract}
\bigskip
\bigskip



\section{Introduction}
Free boundary value problems (BVPs) occur in all branches of applied mathematics and science.
The oldest problem of this type was formulated by Isaac Newton, in book II of his great \lq \lq Principia Mathematica\rq \rq \ of 1687, by considering the optimal nose-cone shape for the motion of a projectile subject to air resistance, see Edwards \cite{Edwards:NNC:1997} or Fazio \cite{Fazio:2014:NTM}.

In the classical numerical treatment of a free BVP a preliminary reduction to a BVP is introduced by considering a new independent variable; see, Stoer and Bulirsch \cite[p. 468]{StoerBBook},
Ascher and Russell \cite{Ascher:1981:RBV}, or Ascher, Mattheij and Russell \cite[p. 471]{AscherBook}.
By rewriting a free BVP as a BVP it becomes evident that the former is always a nonlinear problem; the first to point out the nonlinearity of free BVPs was Landau \cite{Landau:1950:HCM}.
Therefore, in that way free BVPs are BVPs.
In this paper we show that free BVPs invariant with respect to a translation group can be solved non-iteratively by solving a related initial value problem (IVP). 
Therefore in this way those free BVPs are indeed IVPs.
Moreover, we are able to characterize a class of free BVPs that can be solved non-iteratively by solving related IVPs.

The non-iterative numerical solution of BVPs is a subject of past and current research.
Several different strategies are available in literature for the non-iterative solution of BVPs:
superposition \cite[pp. 135-145]{AscherBook}, chasing \cite[pp. 30-51]{Na:1979:CME}, and adjoint operators method \cite[pp. 52-69]{Na:1979:CME} that can be applied only to linear models; parameter differentiation \cite[pp. 233-288]{Na:1979:CME} and invariant imbedding \cite{Meyer:1973:IVM} can be applied also to nonlinear problems. 
In this context transformation methods (TMs) are founded on group invariance theory, see Bluman and Cole \cite{Bluman:1974:SMD}, Dresner \cite{Dresner:1983:SSN}, or Bluman and Kumei \cite{Bluman:1989:SDE}. 
These methods are initial value methods because they achieve the numerical solution of BVPs through the solution of related IVPs.

The first application of a non-iterative TM was given by T\"opfer in \cite{Topfer:1912:BAB} for the Blasius problem, without any consideration of group invariance theory.
T\"opfer's algorithm is quoted in several books on fluid dynamics, see, for instance, Goldstein \cite[pp. 135-136]{Goldstein:1938:MDF}.
Acrivos, Shah and Petersen \cite{Acrivos:1960:MHE} first and Klamkin \cite{Klamkin:1962:TCB} later extended T\"opfer's method respectively to a more general problem and to a class of problems.
Along the lines of the work of Klamkin, for a given problem Na \cite{Na:1967:TBC,Na:1968:FET} showed the relation between the invariance properties, with respect to a linear group of transformation (the scaling group), and the applicability of a non-iterative TM.
Na and Tang \cite{Na:1969:MSH} proposed a non-iterative TM based on the spiral group and applied it to a non-linear heat generation model.
Belford \cite{Belford:1969:IVP} first, and Ames and Adams \cite{Ames:1976:ESE,Ames:1979:NLB} later defined non-iterative TMs for eigenvalue problems. 
A review paper was written by Klamkin \cite{Klamkin:1970:TBV}.
Extensions of non-iterative TM, by requiring the invariance of one and of two or more physical parameters when they are involved in the mathematical model, were respectively proposed by Na \cite{Na:1970:IVM} and by Scott, Rinschler and Na \cite{Scott:1972:FEI}; see also Na \cite[Chapters 8 and 9]{Na:1979:CME}. 
A survey book, written by Na \cite[Chs 7-9]{Na:1979:CME} on the numerical solution of BVP, devoted three chapters to numerical TMs.

As far as  free BVPs are concerned, non-iterative and iterative TMs were proposed by Fazio and Evans \cite{Fazio:1990:SNA}.
Fazio \cite{Fazio:1990:NVT} has shown that we can extend the applicability of non-iterative TMs by rewriting a given free BVP using a variables transformation obtained by linking two different invariant groups. 

However, non-iterative TMs are applicable only to particular classes of BVPs so that they have been considered as \textit{ad hoc} methods, see Meyer \cite[pp. 35-36]{Meyer:1973:IVM}, Na \cite[p. 137]{Na:1979:CME} or Sachdev \cite[p. 218]{Sachdev:1991:NOD}.

The transformation of BVPs to IVPs has also a theoretical relevance.
In fact, existence and uniqueness results can be obtained as a consequence of the invariance properties.
For instance,  for the Blasius problem, a simple existence and uniqueness theorem was given by J. Serrin \cite{Serrin:1970:ETS} as reported by Meyer \cite[pp. 104-105]{Meyer:1971:IMF} or Hastings and McLeod \cite[pp. 151-153]{Hastings:2012:CMO}.
Moreover, using scaling invariance properties the error analysis of the truncated boundary formulation of the Blasius problem was developed by Rubel \cite{Rubel:1955:EET}.
On this topic a first application of a numerical test, defined within group invariance theory, to verify the
existence and uniqueness of the solution of a free BVPs was considered by Fazio in \cite{Fazio:1991:ITM}.
A formal definition of the mentioned numerical test can be found in \cite{Fazio:1997:NTE}.

In this paper we consider the class of free BVPs governed by second order autonomous differential equations, and define, for these problems, a non-iterative TM using the invariance properties of a translation group.
As far as applications of the proposed algorithm are concerned, we solve three problems. 
First we test our method with a problem describing a rope configuration against an obstacle, where we compare the obtained numerical results with the exact solution. 
Then we solve a dynamical problem with a nonlinear force, and a problem related to the optimal length estimate for tubular flow reactors, where in both cases our results are compared to numerical data available in literature. 
Finally, the last section is concerned with concluding remarks pointing out limitations and possible extensions of the proposed approach.

\section{The non-iterative TM}
Let us consider the class of second order free BVPs given by
\begin{align}
&\frac{d^2u}{dx^2}=\Omega\left(u,\frac{du}{dx}\right) \ , \qquad x \in (0, s) \label{gov-t1}\\
&A_1 u(0) +A_2 \frac{du}{dx}(0) = A_3 \ , \label{l-bound-t1}\\
&u(s)=B \quad , \quad \frac{du}{dx}(s)=C\ ,\label{r-bound-t1}
\end{align}
where $A_i$, for $i=1,2,3$, $B$ and $C$ are arbitrary constants, and $s>0$ is an unknown free boundary.
The differential equation (\ref{gov-t1}) and the two free boundary conditions (\ref{r-bound-t1}) are invariant with respect to the translation group
\begin{align}
x^*=x+\mu\quad ; \quad s^*=s+\mu\quad ; \quad u^*=u \ .\label{t1}
\end{align}
By using this invariance, we can define the following non-iterative algorithm for the numerical solution of (\ref{gov-t1})-(\ref{r-bound-t1}):
\begin{itemize}
\item we fix freely a value of $s^*$; 
\item we integrate backwards from $ s^* $ to $ x_{0}^* $ the following auxiliary IVP 
\begin{align}\label{eq:IVP}
&\frac{d^2u^*}{dx^{*2}}=\Omega\left(u^*,\frac{du^*}{dx^*}\right) \nonumber \\[-1ex]
& \\[-1ex]
&u^*(s^*)=B\ , \qquad \frac{du^*}{dx^*}(s^*)=C \ ,	\nonumber
\end{align}
using an {\it event locator} in order to find $ x_{0}^* $ such that 
\begin{equation}\label{eq:evloc} 
A_1 u^*(x_{0}^*) + A_2\frac{du^*}{dx^*}(x_{0}^*) = A_3 \ ;
\end{equation} 
\item finally, through the invariance property,  
we can deduce the similarity parameter
\begin{equation}\label{eq:mu}
\mu = x_{0}^* \ ,
\end{equation} 
from which we get the unknown free boundary
\begin{equation}\label{eq:freeb}
s = s^*-\mu \ . 
\end{equation} 
The missing initial conditions are given by
\begin{equation}\label{eq:mics}
u(0) = u^*(x_{0}^*) \ , \qquad \frac{du}{dx}(0) = \frac{du^*}{dx^*}(x_{0}^*) \ .
\end{equation} 
\end{itemize}

Let us define now a simple event locator suited to the class of problems (\ref{gov-t1})-(\ref{r-bound-t1}).
We consider first the case where
\begin{equation}\label{eq:evloc1}
A_1 u^*(s^*) +A_2 \frac{du^*}{dx^*}(s^*) < A_3 \ .
\end{equation} 
We can integrate the auxiliary IVP (\ref{eq:IVP}) with a constant step size $\Delta x^*$ until
at a given mesh point $x_k^*$ we get
\begin{equation}\label{eq:evloc2}
A_1 u^*(x_k^*) +A_2 \frac{du^*}{dx^*}(x_k^*) > A_3 \ ,
\end{equation} 
and repeat the last step with the smaller step size
\begin{equation}\label{eq:evloc3}
\Delta x_0^* = \Delta x^* \frac{\displaystyle A_3 - A_1 u^*(x_{k-1}^*) -A_2 \frac{du^*}{dx^*}(x_{k-1}^*)}{\displaystyle A_1 u^*(x_k^*) +A_2 \frac{du^*}{dx^*}(x_k^*)- A_1 u^*(x_{k-1}^*) -A_2 \frac{du^*}{dx^*}(x_{k-1}^*)} \ .
\end{equation} 
In defining the last step size in equation (\ref{eq:evloc3}) we use a linear interpolation. 
As a consequence, we have that $x_0^* \approx x_k^* - \Delta x_0^*$.
We notice that the condition imposed by this event locator converges to the correct condition
(\ref{eq:evloc}) as the step size goes to zero, cf. the second column of table~\ref{Fazio:Iacono:tb3}. 

The other case
\begin{equation}\label{eq:evloc4}
A_1 u^*(s^*) +A_2 \frac{du^*}{dx^*}(s^*) > A_3 \ ,
\end{equation} 
can be treated in a similar way.
Of course, also in this second case the last step size is smaller than the previous ones.

In the next section we apply the proposed non-iterative TM to three problems.
The reported numerical results were computed by the classical fourth-order Runge-Kutta's method, reported by Butcher \cite[p. 166]{Butcher}, coupled with the event locator defined above.

\section{The obstacle problem on a string}
\begin{figure}[!hbt]
	\centering
\psfragscanon
\psfrag{u}{$u(x)$}
\psfrag{u0}{$u_0$}
\psfrag{s}{$s$}
\psfrag{du}{$\displaystyle \frac{du}{dx}(x)$}
\psfrag{x}{$x$}
\includegraphics[width=.75\textwidth]{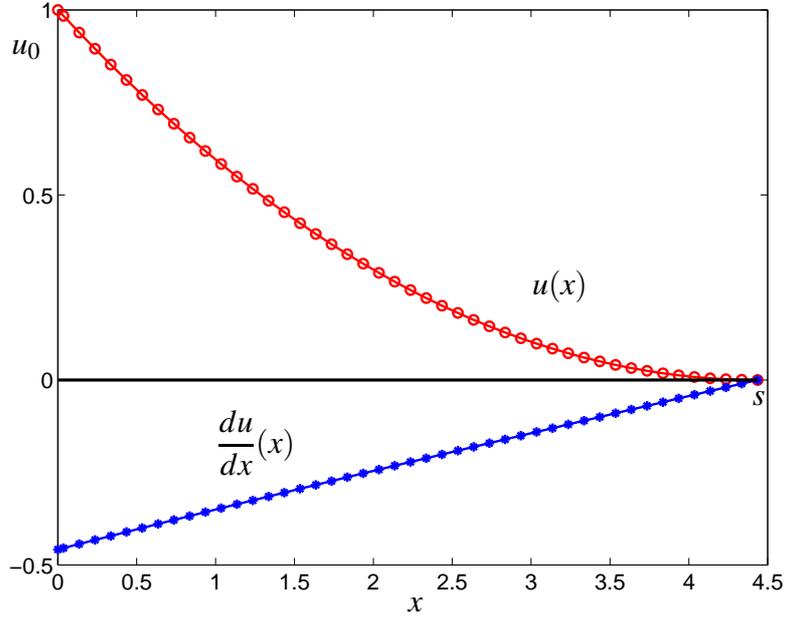}
\caption{Picture of the numerical solution obtained for $\theta=0.1$, $u_0=1$, and $ b = 4.5 $; the obstacle, is superimposed to the $x$ axis, and is displayed by a black solid line.}
\label{fig:string}
\end{figure}

For the obstacle problem on a string, depicted on figure~\ref{fig:string} within the $(x,u)$-plane  where the $ x $ axis is taken overlying to the obstacle, we have to consider the following mathematical model, see Collatz \cite{Collatz:1980:MFB} or Glashoff and Werner, \cite{Glashoff:IMM:1979}
\begin{align}\label{eq:free:string}
&\frac{d^2u}{dx^2}=\theta\left[1+\left(\frac{du}{dx}\right)^2\right]^{1/2} \ , \qquad x \in (0, s) \nonumber \\[-1ex]
&\\[-1ex]
& u(0) = u_0 \ , \qquad u(s)=\frac{du}{dx}(s)=0 \ , \nonumber 
\end{align}
where the positive value of $ \theta $ depends on the string properties.
In this problem we have to find the position of a uniform string of finite length $ L $ under the action of gravity.
The string has fixed end points, say $ (0, u_0) $ and $ (b,0) $, where $ u_0 > 0 $ and $ b > 0 $.
Furthermore, we assume that the condition $ L^2 > \left(u_0^2+b^2\right) $ is fulfilled; this condition allows us to define a free boundary $ s $ for this problem, where $ s $ is the detached rope position from the obstacle.

The free BVP (\ref{eq:free:string}) was solved by the first author in \cite{Fazio:1991:FBT} by iterative methods, namely a shooting method and the iterative extension of the TM derived by using the invariance with respect to a scaling group.

The exact solution of the free BVP (\ref{eq:free:string}) is given by
\begin{eqnarray}\label{eq:exact}
& u(x) = \theta^{-1} \left[\cosh\left(\theta \left(x-s\right)\right) - 1\right] \ , \nonumber \\[-1ex]
& \\[-1ex]
&s = \theta^{-1} \ln \left[\theta u_0+1 + \left(\left(\theta u_0 +1\right)^2 - 1\right)^{1/2}\right] \ , \nonumber
\end{eqnarray}
from this we easily find
\begin{equation}\label{eq:exactdu}
\frac{du}{dx}(0) = \sinh(-\theta s)=\frac{1}{2} \left(e^{-\theta s} - e^{\theta s}\right) \ ,
\end{equation}
and, therefore, for $\theta = 0.1$ and $u_0 = 1$ from equations (\ref{eq:exact})-(\ref{eq:exactdu}) we get the values
\begin{equation}\label{eq:ex}
s=4.356825433 \ , \qquad \frac{du}{dx}(0) = -0.458257569 \ , 
\end{equation}
that are correct to the ninth decimals. 

\begin{table}[!hbt]
\caption{Convergence of numerical results for $\theta= 0.1$.}
\begin{center}{\renewcommand\arraystretch{1.2}
\begin{tabular}{l|cccc}
\hline
{$\Delta x$} & {$\frac{du}{dx}(0)$} & $e_r$ & $s$ & $e_r$ \\
\hline  
$-0.1$  & $-0.458227362$ & $6.59\mbox{D}-05$ & $4.435407932$ & $6.19\mbox{D}-05$ \\
$-0.05$ & $-0.458250809$ & $1.47\mbox{D}-05$ & $4.435621088$ & $1.39\mbox{D}-05$ \\
$-0.025$ & $-0.458255551$ & $4.40\mbox{D}-06$ & $4.435664194$ & $4.14\mbox{D}-06$ \\
$-0.0125$ & $-0.458257313$ & $5.59\mbox{D}-07$ & $4.435680211$ & $5.26\mbox{D}-07$ \\
$-0.00625$ & $-0.458257463$ & $2.31\mbox{D}-07$ & $4.435681576$ & $2.18\mbox{D}-07$ \\
$-0.003125$ & $-0.458257538$ & $6.74\mbox{D}-08$ & $4.435682258$ & $6.43\mbox{D}-08$  \\
$-0.0015625$ & $-0.458257565$ & $8.52\mbox{D}-09$ & $4.435682504$ & $9.05\mbox{D}-09$ \\
\hline
\end{tabular}
}
\end{center}
\label{Fazio:Iacono:tb}
\end{table}

Let us consider a convergence numerical test for our non-iterative TM.
Table~\ref{Fazio:Iacono:tb} reports the obtained numerical results for the missing initial condition and the free boundary value for the free BVP (\ref{eq:free:string}) with $\theta = 0.1$ and $u_0 = 1$, as well as the corresponding relative errors denoted by $e_r$.  
An example of the numerical solutions is shown in figure~\ref{fig:string}.
For this numerical solution, we applied a large step size in order to show the mesh used and to empathize how our event locator reduces the last step.

\section{A dynamical free BVP}
Suppose a particle of unitary mass is moving against a nonlinear force, given by $-1-u-\left(\frac{du}{dx}\right)^2$, from the origin $u=0$ to a final position $u=1$, our goal is to determine the duration of motion $s$ and the initial velocity that assures that the particle is momentarily at rest at $u=1$; see Meyer \cite[pp. 97-99]{Meyer:1973:IVM}.
This problem can be formulated as follows 
\begin{align}\label{eq:free:dyn}
&\frac{d^2u}{dx^2}=-1-u-\left(\frac{du}{dx}\right)^2 \ , \qquad x \in (0, s) \nonumber \\[-1ex]
&\\[-1ex]
& u(0) = 0 \ , \qquad u(s)=1 \ , \qquad \frac{du}{dx}(s)=0 \ , \nonumber 
\end{align}
where $u$ and $x$ are the particle position and the time variable, respectively, on the right hand side of the governing differential equation we have the nonlinear force acting on the particle and $s$ is the free boundary. 

In table~\ref{Fazio:Iacono:tb3} we propose a numerical convergence test for decreasing values of the step size.
\begin{table}[!hbt]
\caption{Convergence of numerical results for the free BVP dynamical model.}
\begin{center}{\renewcommand\arraystretch{1.2}
\begin{tabular}{l|ccc}
\hline
{$\Delta x$} & {$u(0)$} & {$\frac{du}{dx}(0)$} & $s$  \\
\hline  
$-0.1$ & $1.16\mbox{D}-02$ & $3.212263787$ & $0.867662139$ \\
$-0.05$ & $3.54\mbox{D}-03$ & $3.240676696$ & $0.870143219$ \\
$-0.025$ & $4.61\mbox{D}-04$ & $3.2516023692$ & $0.871089372$ \\
$-0.0125$ & $1.90\mbox{D}-04$ & $3.252564659$ & $0.871172452$ \\
$-0.00625$ & $5.42\mbox{D}-05$ & $3.253049203$ & $0.871214290$ \\
$-0.003125$ & $9.25\mbox{D}-06$ & $3.253209165$ & $0.871228100$ \\
$-0.0015625$ & $3.43\mbox{D}-06$ & $3.253229900$ & $0.871229890$ \\
$-0.00078125$ & $5.12\mbox{D}-07$ & $3.253240276$ & $0.871230785$ \\
$-0.000390625$ & $2.01\mbox{D}-07$ & $3.253241381$ & $0.871230881$ \\
$-0.0001953125$ & $4.62\mbox{D}-08$ & $3.253241934$ & $0.871230929$ \\
\hline
\end{tabular}
}
\end{center}
\label{Fazio:Iacono:tb3}
\end{table}

\begin{figure}[!hbt]
\centering
\psfragscanon 
	\psfrag{x}{$x$}
	\psfrag{u}{$u(x)$}
	\psfrag{du}{$\displaystyle \frac{du}{dx}(x)$}
\includegraphics[width=.75\textwidth]{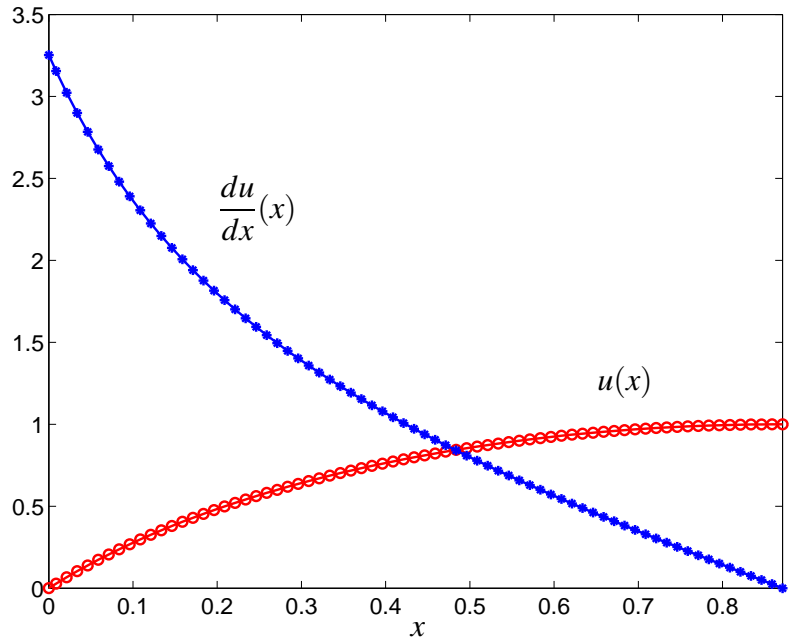} 
	\caption{Numerical solution for the dynamical free BVP (\ref{eq:free:dyn}) obtained with $\Delta x = 0.0125$.}
	\label{fig:free:dyn}
\end{figure}

The obtained results can be contrasted with those reported by Meyer \cite[pp. 97-99]{Meyer:1973:IVM}
where, by using the invariant imbedding method, he found $s=1.2651$ but a value of $u(0) = 0.0163$ instead of $u(0)=0$ as prescribed by the free BVP (\ref{eq:free:dyn}).
The behaviour of the solution can be seen in the figure~\ref{fig:free:dyn}.
Again we applied a large step size in order to show how our event locator reduces the last step.

A free BVP similar to (\ref{eq:free:dyn}) was considered by Na \cite[p. 88]{Na:1979:CME} where the nonlinear force was replaced by $-u \exp(-u)$.
However, in this case it is possible to prove \cite{Fazio:1997:NTE}, using the conservation of energy principle, that the free BVP has countable infinite many solutions, with the missing initial conditions given by
\begin{equation}
\frac{du}{dx}(0) = \pm 0.726967811 \ .
\end{equation}

\section{Length estimation for tubular flow reactors}
\begin{figure}[!hbt]
\centering
\psfragscanon 
	\psfrag{x}{\small$x$}
	\psfrag{A}{\small A}
	\psfrag{B}{\small B}
\framebox{
\includegraphics[width=.75\textwidth]{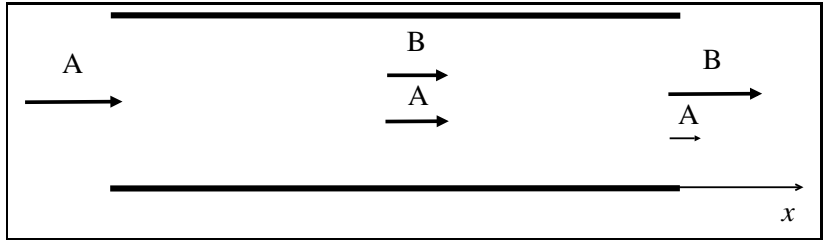} 
}	
\caption{Schematic set-up of a tubular flow reactor.}
	\label{fig:setup}
\end{figure}

Roughly speaking, a tubular flow chemical reactor can be seen as a device where on one side it is introduced a material A that along its passage inside the reactor undergoes a chemical reaction so that at the exit we get a product B plus a residual part of A; see figure~\ref{fig:setup}.
A $n$th order chemical reactor is usually indicated with the notation A$^n \rightarrow$ B.  
A free BVP for a tubular reactor can be formulated as 
\begin{align}\label{eq:free:tubular}
&\frac{d^2u}{dx^2}=N_{Pe}\left(\frac{du}{dx}+Ru^n\right) \ , \nonumber \\[-1ex]
&\\[-1ex]
& u(0)-\frac{1}{N_{Pe}}\frac{du}{dx}(0) = 1 \ , \qquad u(s)=\tau \ , \quad \frac{du}{dx}(s)=0 \ , \nonumber 
\end{align}
where $u(x)$ is the ratio between the concentration of the reactant A at a distance $x$ and the concentration of it at $x=0$, $N_{Pe}$, $R$, $n$ and $\tau$  are, the Peclet group, the reaction rate group, the order of the chemical reaction and the residual fraction of reactant A at exit, respectively. Moreover, $N_{Pe}$ and $R$ are both greater than zero. Finally,  for the free BVP (\ref{eq:free:tubular}), the free boundary $s$ is the length of the flow reactor we are trying to estimate.

This is an engineering problem that consists in determining the optimal length of a tubular flow chemical reactor with axial missing and has been already treated by Fazio in \cite{Fazio:1991:ITM}, through an iterative TM, whereas Fazio in \cite{Fazio:1992:NLE} made a comparison between the results obtained with a shooting method and the upper bound of the free boundary value obtained by a non-iterative TM. 

Here, for the sake of comparing the numerical results, we fix the parameters as follows: $N_{Pe}=6$, $R=2$, $n=2$, and $\tau=0.1$.
We apply the algorithm outlined above to the numerical solution of the free BVP (\ref{eq:free:tubular}).
Table~\ref{Fazio:Iacono:tb1} shows a numerical convergence test for decreasing values of the step size.
\begin{table}[!hbt]
\caption{Convergence of numerical results. $N_{pe}=6$, $R=2$, $n=2$, and $\tau= 0.1$.}
\begin{center}{\renewcommand\arraystretch{1.2}
\begin{tabular}{l|ccc}
\hline
{$\Delta x$} & {$u(0)$} & {$\frac{du}{dx}(0)$} & $s$  \\
\hline  
$-0.1$ & $0.829314641$ & $-1.008175212$ & $5.117905669$ \\
$-0.05$ & $0.830537187$ & $-1.010745699$ & $5.119104349$ \\
$-0.025$ & $0.831147822$ & $-1.012077034$ & $5.119707352$ \\
$-0.0125$ & $0.831227636$ & $-1.012251496$ & $5.119786158$ \\
$-0.00625$ & $0.831267467$ & $-1.012338738$ & $5.119825502$ \\
$-0.003125$ & $0.831271635$ & $-1.012347868$ & $5.119829619$ \\
$-0.0015625$ & $0.831273719$ & $-1.012352436$ & $5.119831678$ \\
$-0.00078125$ & $0.831274182$ & $-1.012353449$ & $5.119832135$ \\
$-0.000390625$ & $0.831274327$ & $-1.012353767$ & $5.119832278$ \\
$-0.0001953125$ & $0.831274348$ & $-1.012353814$ & $5.119832299$ \\
\hline
\end{tabular}
}
\end{center}
\label{Fazio:Iacono:tb1}
\end{table}

The obtained numerical results are reported on table~\ref{Fazio:Iacono:tb2} and compared with numerical results available in literature.
\begin{table}[!hbt]
\caption{Comparison of numerical results for the tubular flow reactor model.}
\begin{center}{\renewcommand\arraystretch{1.2}
\begin{tabular}{lccc}
\hline
& {$u(0)$} & {$\frac{du}{dx}(0)$} & {$s$} \\  
\hline
iterative TM \cite{Fazio:1991:ITM} & $0.831280$ & $-1.012298$ & $5.121648$ \\
shooting method \cite{Fazio:1992:NLE} & $0.831274$ & $-1.012354$ & $5.119832$ \\
non-iterative TM & $0.831274$ & $-1.012354$ & $5.119832$ \\
\hline
\end{tabular}
}
\end{center}
\label{Fazio:Iacono:tb2}
\end{table}

\begin{figure}[!hbt]
\centering
\psfragscanon 
	\psfrag{x}{$x$}
	\psfrag{u}{$u(x)$}
	\psfrag{du}{$\displaystyle \frac{du}{dx}(x)$}
\includegraphics[width=.75\textwidth]{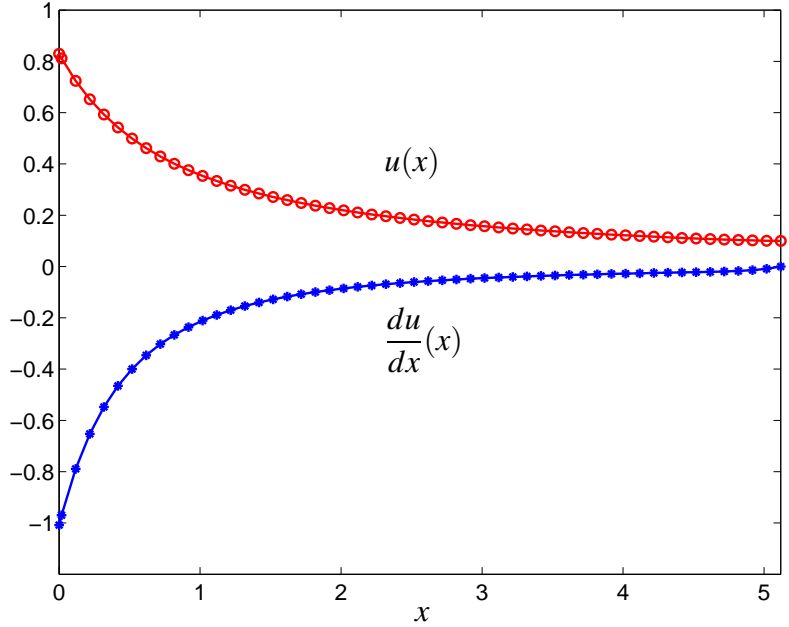} 
	\caption{Numerical solution for length estimation of a tubular flow reactor obtained with $\Delta x = 0.1$.}
	\label{fig:tubular}
\end{figure}

As it is easily seen the computed values are in good agreement with the ones found in \cite{Fazio:1991:ITM} and \cite{Fazio:1992:NLE}. 
The behaviour of the solution can be seen in figure~\ref{fig:tubular}.
Once again, we used a large step size to make clear how our event locator reduces the last step size.

\section{Conclusion}
In closing, we can outline some further implications coming out from this work. 
First of all, the algorithm proposed in this paper can be extended to free BVPs governed by a system of first order autonomous differential equations belonging to the general class of problems
\begin{align}
& {\displaystyle \frac{d{\bf u}}{dx}} = {\bf q} \left({\bf u}\right)
\ , \quad x \in [0, \infty) \ , \nonumber \\[-1.5ex]
\label{psystem} \\[-1.5ex]
& u_j(0) = u_{j0} \ ,  \qquad {\bf u} (s) = {\bf u}_s
\ ,  \nonumber
\end{align}
where ${\bf u}(x), {\bf u}_s \in \RR^d$, $ {\bf q}: \RR^d \rightarrow~\RR^d $, with $d \ge 1$,
$j \in \{1, \dots , d \}$, $ u_{j0} $ and all components of ${\bf u}_s $ are given constants and $s$ is the free boundary.

Moreover, our algorithm can be applied by using an integrator from the MATLAB ODE suite written by Samphine and Reichelt \cite{Shampine:1997:MOS}, and available with the latest releases of MATLAB, with the {\it event locator} option command set in 
\[
\mbox{options = odeset('Events',@name)}
\]
where \lq \lq name\rq \rq \ is an external, problem dependent, {\it event function}. 

As mentioned in the introduction, the first application of a non-iterative TM was defined by T\"opfer in \cite{Topfer:1912:BAB} more than a century ago.
In this paper, by considering the invariance with respect to a translation group, we have investigated a possible way to solve a large class of free BVPs by a non-iterative TM.

However, it is a simple matter to show a differential equation not admitting any group of transformations: e.g. the differential equation considered by Bianchi \cite[pp. 470-475]{Bianchi:1918:LTG}.
Consequently, it is easy to realize that non-iterative TMs cannot be extended to every BVPs.
Therefore, non-iterative TMs are \textit{ad hoc} methods. 
Their applicability depends on the invariance properties of the governing differential equation and the given boundary conditions.

On the other hand, free BVPs governed by the most general second order differential equation, in normal form, can be solved iteratively by extending a scaling group via the introduction of a numerical parameter so as to recover the original problem as the introduced parameter goes to one, see Fazio \cite{Fazio:1990:SNA,Fazio:1990:NVT,Fazio:1997:NTE,Fazio:1998:SAN}.
The extension of this iterative TM to problems in boundary layer theory has been considered in \cite{Fazio:1994:FSEb,Fazio:1996:NAN,Fazio:2009:NTM,Fazio:2013:BPF}.
\ Moreover, a further extension to the sequence of free BVPs obtained by a semi-discretization of parabolic moving boundary problems was repoted in \cite{Fazio:2001:ITM}.

\end{document}